\documentclass[12pt]{article}

\usepackage{amsmath}
\usepackage{amssymb}

\usepackage{a4wide}

 \newcommand{\bs}[1]{\boldsymbol{#1}}
 \newcommand{\floor}[1]{\lfloor{#1}\rfloor}
 \newcommand{\ceiling}[1]{\lceil{#1}\rceil}

 \newcommand{\RR}{\mathbb R}
 
 \newcommand{\ZZ}{\mathbb Z}

 \newcommand{\cM}{\mathcal M}
 \newcommand{\cP}{\mathcal P}
 \newcommand{\cB}{\mathcal B}
 \newcommand{\cR}{\mathcal R}

 \newcommand{\es}{\varnothing}

 \newcommand{\s}{{\scriptscriptstyle <}}
 \newcommand{\g}{{\scriptscriptstyle >}}

 \newtheorem{theorem}{Theorem}
 \newtheorem{lemma}{Lemma}
 \newtheorem{prop}{Proposition}
 \newtheorem{coro}{Corollary}
 \newtheorem{fact}{Fact}

 \newcommand{\qed}{\hfill $\square$}
 \newcommand{\btimes}{\mbox{\huge \raisebox{-0.2ex}{$\times$}}}

\begin{document}

\centerline{\Large \bf Recombination semigroups on measure spaces}

\vspace{15mm}
        
\centerline{\sc Michael Baake}
\vspace{5mm} 

{\small
\begin{center}
    Fakult\"at f\"ur Mathematik, Univ.\ Bielefeld, \\
    Box 10 01 31, 33501 Bielefeld, Germany  
\end{center}
}

\vspace{20mm}
\begin{abstract} 
  The dynamics of recombination in genetics leads to an interesting
  nonlinear differential equation, which has a natural generalization
  to a measure valued version. The latter can be solved explicitly
  under rather general circumstances. It admits a closed formula for
  the semigroup of nonlinear positive operators that emerges from the
  forward flow and is, in general, embedded in a multi-parameter
  semigroup.

\end{abstract}

\vspace{30mm}
\hphantom{Key Words:}           
measure-valued dynamical systems, M\"obius inversion

\smallskip
\noindent
2000 MCS: $\;$92D10, 34L30 (primary); 37N30, 06A07, 60J25 (secondary)

\clearpage

\section{Introduction}

The deterministic limit of the stochastic process of recombination
in population genetics leads to an interesting nonlinear ODE system
that has been studied for a long time. The first major advances to
understand the classical system are due to H.\ Geiringer
\cite{Gei} in the 1940s (though her formulation was slightly
different). A full characterization of solutions (though not
in explicit form) was later obtained by Lyubich, compare \cite{Lyu}
and references given there, and \cite{Reinhard} for general
background material.

Motivated by this problem in genetics, a relevant subclass was
considered in \cite{EB}, where an explicit solution to the resulting 
large system of nonlinear differential equations was constructed.
Soon after, this was reformulated as a measure valued
differential equation and solved for a more general class of state
spaces \cite{reco}. In this context, the main focus was not on the
functional analytic properties, but on the solution and the
interaction with other genetic processes, in particular mutation and
selection.

On the other hand, nonlinear semigroups are rarely known explicitly,
and if so, this usually rests upon the transformability of the system
to a linear one. With this observation in mind, it seems worthwhile to
also investigate the semigroup aspect separately. To do so, we will
first consider the most elementary semigroup constituents separately.
Among them, we will identify mutually commuting ones, which are then
used to reconstruct the solution to the full recombination equation of
\cite{reco} in a simpler than the previous way. Furthermore,
generalizations are possible that are less relevant in biological
applications, but illustrate that the mathematical analysis may still
be pushed considerably.

In view of the origin of the problem, the formulation is oriented
towards the application in mathematical biology. Consequently, the
formulation is slightly more expository than necessary for a purely
mathematical audience. In a separate section, an abstract
reformulation of the central observation is added, together with a
possible generalization to a wider class of nonlinear operators.

\section{Mathematical setting}

If $X$ is a locally compact space (by which we always mean to include
the Hausdorff property), we denote by $\cM(X)$ the Banach space of
{\em finite\/} regular Borel measures, equipped with the usual
variation norm $\|.\|$, i.e., $\|\omega\| := |\omega| (X)$ where
$|\omega|$ denotes the total variation measure of $\omega$.  $\cM(X)$
is a Banach lattice, and we use $\cM_+(X)$ (resp.\ $\cP(X)$) to denote
the closed convex cone of positive measures (resp.\ the closed simplex
of probability measures).  For positive measures $\omega$, one has
$\|\omega\| = \omega(X)$. We will mainly be concerned with the
situation that $X$ is itself a product space, e.g., $X=X_1 \times X_2$
with $X_i$ locally compact. Let us first recall a well-known result,
compare \cite{Berb} and \cite[Fact~1]{reco} for details.
\begin{fact} \label{baire}
  Let $\nu,\nu'$ be two regular Borel measures on the locally compact
  product space $X=X_1\times X_2$, and let $\nu,\nu'$ coincide on all
  ``rectangles'' $E_1\times E_2$ where $E_1$ and $E_2$ run through the
  Borel sets of $X_1$ and $X_2$, respectively. Then, $\nu=\nu'$, i.e.,
  $\nu(E)=\nu'(E)$ for all Borel sets\/ $E$ of\/ $X$.  \qed
\end{fact}

The state space $X$ that we need will have a product structure,
described on the basis of sites or {\em nodes}. For later convenience,
we use $N=\{0,1,\dots ,n\}$ for the set of nodes, i.e., we start
counting with $0$ here. To node $i$, we attach the locally compact
space $X_i$, and our state space is then
\begin{equation} \label{seq-space}
   X \; = \; X_0 \times X_1 \times \ldots \times X_n
\end{equation}
which is still locally compact. Our dynamics will evolve in the Banach
space $\cM(X)$ which canonically contains the tensor product
$\cM(X_0)\otimes \ldots \otimes \cM(X_n)$ and also its completion, the
corresponding projective tensor product. In fact, product measures
$\omega = \omega^{}_0 \otimes \ldots \otimes \omega^{}_n$ with
$\omega_i \in \cM(X_i)$ will play an important role below.

The main reason for using the above set $N$ of nodes is that we will
need {\em ordered partitions\/} of $N$ which are uniquely specified by
a set of cuts or {\em crossovers\/}. The possible cut positions are at
the {\em links\/} between nodes, which are denoted by half-integers,
i.e., by elements of the set $L=\{\frac{1}{2}, \frac{3}{2}, \dots
,\frac{2n-1}{2}\}$.  We use Latin indices for nodes and Greek indices
for links, and the implicit rule will be that $\alpha=\frac{2i+1}{2}$
is the link between nodes $i$ and $i+1$.

With this notation, the ordered partitions of $N$ are in one-to-one
correspondence with the subsets of $L$ as follows. If
$A=\{\alpha^{}_1,\dots ,\alpha^{}_p\}\subset L$, with $\alpha_i <
\alpha_{i+1}$, let $N^{}_A = \{A^{}_0, A^{}_1, \dots , A^{}_p\}$
denote the ordered partition
$$A^{}_0 = \{0,\dots ,\floor{\alpha^{}_1}\} \, , \, A^{}_1 =
\{\ceiling{\alpha^{}_1},\dots ,\floor{\alpha^{}_2}\}\, , \; \dots \; ,
A^{}_p = \{\ceiling{\alpha^{}_{p}},\dots , n \} $$
where
$\floor{\alpha}$ ($\ceiling{\alpha}$) is the largest integer below
$\alpha$ (the smallest above $\alpha$).  So, while $A\subset L$, one
has $A_i\subset N$ for $0\le i\le p$.  In particular, we have
$N^{}_{\es} = N$ and $N^{}_L = \{\{0\},\dots,\{n\}\}$.  With this
definition, it is clear that $N^{}_B$ is a {\em refinement\/} of
$N^{}_A$ if and only if $A\subset B$. Consequently, the lattice of
ordered partitions of $N$ now corresponds to the Boolean algebra of
all subsets of the finite set $L$, denoted by $\cB(L)$, cf.\ 
\cite[Ch.\ I.2]{Aigner}. We prefer this notation to that with
partitions, as it is easier to deal with.  If $A\subset B$, we will
write $B\!-\!A$ for $B\setminus A$, and $\overset{\;\_\!\_}{A}$ for
the set $L\!-\!A$.

Let $X,Y$ be two locally compact spaces with attached measure spaces
$\cM(X)$ and $\cM(Y)$. If $f\!: X\to Y$ is a continuous function and
$\omega\in\cM(X)$, $f.\omega := \omega\circ f^{-1}$ is an element of
$\cM(Y)$, where $f^{-1}(y) := \{x\in X\mid f(x)=y\}$ means the
preimage of $y\in Y$ in $X$, with obvious extension to $f^{-1}(B)$,
the preimage of a subset $B\subset Y$ in $X$.  Due to the continuity
of $f$, $f^{-1}(B)$ is a Borel set in $X$ if $B$ is a Borel set in
$Y$.

{}From now on, let $X=X_0\times\ldots\times X_n$ and let $N$ and $L$
always denote the set of nodes and links as introduced above. Let
$\pi^{}_i \! : X \to X_i$ be the canonical {\em projection}, which is
continuous. It induces a mapping from $\cM(X)$ to $\cM(X_i)$ by
$\omega\mapsto\pi^{}_i .\omega$, where $(\pi^{}_i .\omega) (E) =
\omega(\pi^{-1}_i (E))$, for any Borel set $E\subset X_i$. By (slight)
abuse of notation, we will use the symbol $\pi^{}_i$ also for this
induced mapping.  It is clear that $\pi_i$ is linear and preserves the
norm of positive measures. In particular, it maps $\cP(X)$ to
$\cP(X_i)$ and may then be understood as marginalization.

Likewise, for any index set $I\subset N$, one defines a projector
$\pi^{}_I\!: \cM(X) \to \cM(X_I)$ with $X_I:=\btimes_{\!i\in I}\,X_i$.
With this notation, $X_N=X$. We will frequently also use the
abbreviation $\pi^{}_{\s\alpha}$ for the projector $\pi^{}_{\{1,...
  ,\floor{\alpha}\}}$, and $\pi^{}_{\g\alpha}$ for
$\pi^{}_{\{\ceiling{\alpha},... ,n\}}$.

This enables us to introduce a class of nonlinear operators, called
{\em recombinators\/} from now on. For $\alpha\in L$, we first define
an elementary recombinator $R_{\alpha}\! : \, \cM(X)\to\cM(X)$ by
$R_{\alpha}(0)=0$ and, if $\omega \neq 0$, by
\begin{equation} \label{reco-op1}
   R_{\alpha} (\omega) \; := \; \frac{1}{\|\omega\|}\,
     \bigl((\pi^{}_{\s\alpha}.\omega) \otimes 
           (\pi^{}_{\g\alpha}.\omega)\bigr)
\end{equation}
which is a (partial) product measure. Here and in what follows, we
tacitly identify (if necessary) a product measure with its unique
extension to a regular Borel measure on $X$, which is justified by
Fact \ref{baire}.  These elementary recombinators satisfy a number of
useful properties \cite{reco}. In particular, $\|R_{\alpha}(\omega)\|
\le \|\omega\|$ for all $\omega\in\cM(X)$, $R_{\alpha}$ is globally
Lipschitz on $\cM(X)$ with Lipschitz constant $\le 3$, and
$R_{\alpha}(a\,\omega) = |a| R_{\alpha}(\omega)$ for all $a\in\RR$ and
$\omega\in\cM(X)$.  Moreover, $R_{\alpha}$ maps $\cM_+(X)$ into itself
and preserves the norm of positive measures, so that $\cP(X)$ is
mapped into itself, too. Finally, when restricted to $\cM_+(X)$, the
elementary recombinators are idempotents and commute, i.e.,
$R^2_{\alpha} = R^{}_{\alpha}$ and $R_{\alpha} R_{\beta} = R_{\beta}
R_{\alpha}$ for all $\alpha,\beta\in L$.

As the next step, let $A=\{\alpha^{}_1, \dots, \alpha^{}_p\} \subset
L$, with $\alpha_i < \alpha_{i+1}$, and let $N_A$ be the corresponding
ordered partition of $N$, as explained above. Then, we define the
(composite) recombinators $R_A$ by
\begin{equation} \label{composite}
   R^{}_{A}(\omega) \; := \; \frac{1}{\|\omega\|^p} \;
            \bigotimes_{i=0}^{p} \big( \pi^{}_{A_i}.\omega\big),
\end{equation}
again with the continuous extension $R^{}_A (0) = 0$.  Note that
$R^{}_{\es}=\bs{1}$ and $R^{}_{\{\alpha\}} = R^{}_{\alpha}$ in this
notation. As a direct consequence of the definition, or by a simple
induction argument based on the properties of the elementary
recombinators, one obtains the following result, compare \cite{reco}
for details.
\begin{prop} \label{alg-rules}
   Let $R^{}_A$ be the $($composite\/$)$ recombinator attached to a subset 
   $A$ of\/ $L$. Then, the following assertions hold.
\begin{itemize}
\item[$1.$] $\|R^{}_A(\omega)\| \le \|\omega\|$ for all $\omega\in\cM(X)$,
   and $R^{}_A$ is globally Lipschitz on $\cM(X)$. 
\item[$2.$] $R^{}_A$ is positive homogeneous of degree $1$, i.e., 
   $R^{}_A(a\,\omega) = |a| R^{}_A(\omega)$ for all $a\in\RR$
   and $\omega\in\cM(X)$.
\item[$3.$] $R^{}_A$ maps $\cM_+(X)$ into itself and preserves the norm of 
   positive measures. In particular, it maps\/ $\cP(X)$ into itself.
 \item[$4.$] On $\cM_+(X)$, one has $R^{}_A = \prod_{\alpha\in A}
   R^{}_{\alpha}$, and the recombinators satisfy the equation
   $$ R^{}_{G} R^{}_{H} \; = \; R^{}_{G\cup H} \, , $$ 
   for arbitrary $G,H\subset L$.   \qed
\end{itemize}
\end{prop}
Note that Part 4, when applied to singleton sets, comprises both the
idempotency and the commutativity of the elementary recombinators.

\section{Recombination dynamics}

Originally motivated by a problem in biology, compare \cite{EB}, we
are interested in the solutions of the nonlinear ODE
\begin{equation}  \label{gen-reco-eq}
   \dot{\omega} \; = \;  \varPhi(\omega) \; := \; \sum_{G\subset L}
   \varrho^{}_{G}\, \big(R^{}_{G} - 1\big) (\omega)
\end{equation}
on the Banach space $\cM(X)$, with the restriction $\varrho^{}_{G}\ge
0$.  In general, even though the Picard-Lindel\"of theorem applies
($\varPhi$ is globally Lipschitz) and the time evolution is given by
the flow of the ODE, it does not seem possible to write down a closed
formula for the corresponding semigroup in forward time.

However, if $\varrho^{}_{G} > 0$ only for special $G$, namely the
singleton sets $\{\alpha\}$ with $\alpha\in L$, it was shown in
\cite{reco} that the solution in forward time can be given in closed
form, if the initial condition is a positive measure.  This amounts to
constructing an explicit formula for the nonlinear flow of this ODE on
the positive cone $\cM_+(X)$. In what follows, we will give an
independent derivation of this, and an extension to further solvable
cases, independent of its biological relevance. Let us start with a
general property.
\begin{prop} \label{reco-dgl}
  The abstract Cauchy problem of the ODE $(\ref{gen-reco-eq})$ has a
  unique solution. Furthermore, $\cM_+(X)$ is positive invariant under
  the flow, with the norm of positive measures preserved.  In
  particular, $\cP(X)$ is positive invariant.
\end{prop}
{\sc Proof}: Existence and uniqueness of the solution follows from the
Picard-Lindel\"of theorem on Banach spaces, see \cite[Thm.~7.6 and
Remark~7.10]{Amann}, because $\varPhi$ is globally Lipschitz as a
consequence of Part 1 of Proposition~\ref{alg-rules}.

Let $\nu\in\cM_+(X)$, i.e., $\nu(E)\ge 0$ for all Borel sets $E\subset
X$.  Let $E$ be any Borel subset of $X$ such that $\nu(E)=0$. Then
$$
\varPhi(\nu)(E) \; = \; \sum_{G\subset L} \varrho^{}_{G}
R^{}_{G}(\nu)(E) \; \ge \; 0 $$
because each $R^{}_{G}(\nu)$ is a
positive measure and all $\varrho^{}_{G} \ge 0$ by assumption, so
$\varPhi$ satisfies the positive minimum principle, compare
\cite{Arendt}. By a continuity argument, see p.\ 235, Thm.~16.5 and
Remark 16.6 of \cite{Amann}, we obtain the invariance of the cone
$\cM_+(X)$ under the flow in forward time.

{}Finally, observe that $\varPhi(\nu) (X) = 0$ for all
$\nu\in\cM_+(X)$, which is a consequence of Part 3 of
Proposition~\ref{alg-rules}.  Let $\omega^{}_0 \in\cM_+(X)$ be the
initial condition, and $\omega^{}_t$ the corresponding solution. So,
$\omega^{}_t \in\cM_+(X)$ for all $t\ge 0$ by the previous argument,
hence $\|\omega^{}_t\| = \omega^{}_t (X)$.  This implies $\frac{\rm
  d}{{\rm d}t} \|\omega^{}_t\| = \varPhi (\omega^{}_t) (X) = 0$, so
$\|\omega^{}_t\|\equiv \|\omega^{}_0\|$ for all $t\ge 0$.  \qed

\subsection{Time evolution with one recombinator}

Proposition~\ref{reco-dgl} is important because it allows us to
concentrate on the positive cone, and to study the forward flow, which
is a nonlinear positive semigroup.  Our first step is to analyze the
explicit solution, in forward time, of the ODE
\begin{equation} \label{single-rec}
   \dot{\omega} \; = \; \varrho^{}_{A} \big( R^{}_{A} - 1 \big) (\omega)
\end{equation}
for an arbitrary, but fixed $A\subset L$, with $\varrho^{}_A > 0$.
\begin{lemma}  \label{help1}
    Let $\nu\in\cM_+(X)$ and $a\in [0,1]$. Then, we have
    $$
    R^{}_{A} \big(a \nu + (1-a) R^{}_{A}(\nu) \big) \; = \;
    R^{}_{A} (\nu)\, . $$
    So, $R^{}_A$ acts linearly on this type of
    convex combination, and we also get
    $$\big( R^{}_{A} - 1 \big) \big(a \nu + (1-a) R^{}_{A}(\nu) \big)
      \; = \; a \big( R^{}_{A} -1 \big) (\nu)\, . $$
\end{lemma}
{\sc Proof}: Since $R^{}_{A}$ is positive homogeneous of degree 1, it
is sufficient to prove the claim for $\nu\in\cP(X)$, where
$\|\nu\|=1$.  If $A=\{\alpha^{}_1,\dots, \alpha^{}_p \} \subset L$,
and $N^{}_{A} = \{A^{}_0, A^{}_1, \dots , A^{}_p\}$, with $A_i\subset
N$, is the corresponding partition of $N$, we obtain
\begin{eqnarray*}
   \lefteqn{\hspace*{-5mm}R^{}_{A} (a \nu + (1-a) R^{}_{A}(\nu) )}\\
   & = & \bigotimes_{i=0}^{p} 
   \big( \pi^{}_{A_i} . (a \nu + (1-a) R^{}_{A}(\nu) ) \big) 
   \; = \; \bigotimes_{i=0}^{p} \big( a \, \pi^{}_{A_i} . \nu +
         (1-a)\, \pi^{}_{A_i} . R^{}_{A}(\nu) \big) \\
   & = & \bigotimes_{i=0}^{p} \big( a \, \pi^{}_{A_i} . \nu +
         (1-a)\, \pi^{}_{A_i} . \nu \big) 
   \; = \; \bigotimes_{i=0}^{p}\, (\pi^{}_{A_i} . \nu) \,\; = \,\;
         R^{}_{A} (\nu)\, .
\end{eqnarray*}
This proves the first claim, while a verification of the second
is now straightforward.  \qed

\begin{prop}  \label{single-rec-group}
  Let the initial condition $\omega^{}_0$ for equation
  $(\ref{single-rec})$ be a positive measure. Then, the solution for
  $t\ge 0$ is $ \omega^{}_t = \varphi^{A}_{t} (\omega^{}_0)$ where
  $\{\varphi^{A}_{t} \mid t\ge 0 \}$ is a semigroup of nonlinear
  positive operators. They are explicitly given by
  $$
  \varphi^{A}_{t} \; = \; e^{-\varrho^{}_{A} t}\, \bs{1} + \big(1 -
  e^{-\varrho^{}_{A} t} \big) R^{}_{A} \; = \; R^{}_A -
  e^{-\varrho^{}_{A} t}\, \big(R^{}_A - \bs{1}\big) $$
  and, on
  $\cM_+(X)$, satisfy the equation
$$ \frac{\rm d}{{\rm d}t}\, \varphi^{A}_{t} \; = \;
   \varrho^{}_{A} (R^{}_{A} - \bs{1}) \circ \varphi^{A}_{t} \, . $$
Finally, $\omega^{}_t \stackrel{\|.\|\,}{\longrightarrow}R^{}_A(\omega^{}_0)$
as $t\to\infty$, for all $\omega^{}_0\in\cM_+(X)$, with the deviation from
the limit decaying exponentially fast in time.
\end{prop}
{\sc Proof}: Let $\omega^{}_0$ be a positive measure, and
$\omega^{}_t = e^{-\varrho^{}_{A} t}\,\omega^{}_0 +
 \big(1 - e^{-\varrho^{}_{A} t} \big) R^{}_{A} (\omega^{}_0)$.
Then, we have
$$   \dot{\omega}^{}_t \; = \; \varrho^{}_{A}\, e^{-\varrho^{}_{A} t}
     \big(R^{}_{A} - 1\big) (\omega^{}_0)\, .  $$
On the other hand, since $e^{-\varrho^{}_{A} t} \in [0,1]$, we can
apply Lemma~\ref{help1} to verify that also
$$   \varrho^{}_{A} \big(R^{}_{A} - 1\big) (\omega^{}_t)  \; = \;
     \varrho^{}_{A} e^{-\varrho^{}_{A} t}
     \big(R^{}_{A} - 1\big) (\omega^{}_0) $$
which establishes that $\omega^{}_t$ is indeed a solution of the
Cauchy problem with initial condition $\omega^{}_0$, while 
uniqueness follows from Proposition~\ref{reco-dgl}.

Let $\{\varphi^{A}_{t}\mid t\ge 0\}$ denote the corresponding forward
flow, where $\varphi^{A}_{t}$ then is the nonlinear operator stated in
the proposition. It maps $\cM_+(X)$ into itself, preserves the norm of
positive measures, and satisfies the equation
$$ \varphi^{A}_{t} \circ \varphi^{A}_{s} \; = \; \varphi^{A}_{t+s} $$
by the general properties of the flow. Since the formula
for the solution is valid for all initial conditions $\omega^{}_0 \in
\cM_+(X)$, the formula for the derivative is correct on $\cM_+(X)$.  

The statement about the norm convergence of $\omega^{}_t$ to the
product measure $R^{}_A(\omega^{}_0)$ follows from the form of
$\varphi^{}_{t}$ by standard arguments.  \qed

\smallskip It is somewhat surprising that such a simple formula for
the solution emerges. It admits a probabilistic interpretation as
follows. If $\varrho^{}_A$ is the rate of the recombination process
with $R^{}_A$, the term $e^{-\varrho^{}_{A} t}$ (resp.\ $(1 -
e^{-\varrho^{}_{A} t})$) is the probability that recombination has not
yet taken place (resp.\ has happened at least once) until time $t$.
This can be substantiated by viewing the ODE (\ref{single-rec}) as the
deterministic limit of an underlying stochastic process, e.g., in the
spirit of \cite[Sec.~11]{EK}.

Let us also note that the formula for $\varphi^{A}_{t}$ would formally
emerge from expanding $\exp\big(\varrho^{}_A t (R^{}_A - \bf{1})\big)$
as if $R^{}_A$ were a linear idempotent, which it isn't.
Nevertheless, the result is still correct as a consequence of
Lemma~\ref{help1}, because linearity on those special convex
combinations is all that is needed to derive the formula properly.

Still, the right hand side of (\ref{single-rec}) is genuinely
nonlinear. In general, it is well known, compare \cite[p.~91]{Engel},
that $T^A_t (f) := f\circ \varphi^A_t$ defines a semigroup of linear
operators on the dual space of $\cM(X)$. One can also define the
corresponding generator, but this approach does not seem to help in
understanding the result of Proposition~\ref{single-rec-group}.

However, assume $\omega^{}_0\in\cM_+(X)$ and define the two signed
measures $\nu^{}_1(t) = R^{}_A (\omega^{}_t)$ and $\nu^{}_2(t) =
\big(R^{}_A - 1\big) (\omega^{}_t)$. Then, one has $\omega^{}_t =
\nu^{}_1(t) - \nu^{}_2(t)$, but also $\nu^{}_1(t)\equiv R^{}_A
\omega^{}_0$ and $\nu^{}_2(t) = e^{-\varrho^{}_{A} t}\, \big(R^{}_A -
1\big) (\omega^{}_0)$, by an application of Lemma~\ref{help1}.
Consequently, these measures satisfy the differential equations
$$ \dot{\nu}^{}_1 \equiv 0 \qquad \mbox{and} \qquad
   \dot{\nu}^{}_2 = -\varrho^{}_A\, \nu^{}_2  $$
which shows that the ODE (\ref{single-rec}), when restricted to
the cone $\cM_+(X)$, is equivalent to a system of two linear ODEs
together with special initial conditions.
Also, $\nu^{}_1$ is the equilibrium, and $\|\nu^{}_2 (t)\| \to 0$
exponentially fast as $t\to\infty$.

\subsection{Abstract reformulation and possible extension}

In view of other applications, it might be worth to extract the
underlying structure of Proposition~\ref{single-rec-group} as follows.
\begin{theorem} \label{idea}
   Let $K$ be a closed convex subset of a Banach space $B$, and
   let $\cR\!:\, K\longrightarrow K$ be a $($nonlinear\/$)$ Lipschitz
   map which satisfies
\begin{equation} \label{gen-cond}
   \cR \big(a x + (1-a) \cR(x)\big) \; = \; \cR(x)
\end{equation}
   for all $a\in [0,1]$ and all $x\in K$.
   Let $\varrho\ge 0$ be arbitrary.

   Then, the $($nonlinear\/$)$ Cauchy problem
\begin{equation} \label{abstract-Cauchy}
   \dot{x} \; = \; \varrho \big(\cR - \bs{1}\big)(x) \; , \quad
   x(0) = x^{}_0 \in K \, ,
\end{equation}
   has the unique solution
   $x(t) = e^{-\varrho t} x^{}_0 + (1 - e^{-\varrho t}) \cR(x^{}_0)$
   for $t\ge 0$, and the entire forward orbit remains in $K$.
\end{theorem}
\noindent {\sc Proof}:
As above, the uniqueness follows from Lipschitz continuity.  The
formula for the solution is again verified by direct differentiation,
using \eqref{gen-cond}. Note that this assumed property comprises the
linearity on special convex combinations and the idempotency of $\cR$
(by setting $a = 0$).  Since $K$ is a closed convex set which is
mapped into itself by $\cR$, its forward invariance under the flow
follows by standard arguments, compare \cite{Amann}. \qed

\bigskip Several generalizations seem possible at first sight, but
most of them break down rather quickly as soon as genuine nonlinearity
of $\cR$ sets in. First, consider the ansatz
\begin{equation} \label{extend-cond}
    \cR \big(a x + (1-a) \cR^n(x)\big) \; = \; \cR(x) ,
\end{equation}
for some $n > 1$, and again for all $a\in [0,1]$ and all $x\in K$.
This includes the relation $\cR^{n+1}(x) = \cR(x)$ for $x\in K$, a
generalization of the idempotency used previously.  Here, an explicit
solution of the Cauchy problem \eqref{abstract-Cauchy} would contain
all the terms $x^{}_{0}, \cR(x^{}_0), \ldots , \cR^n(x^{}_0)$, so that
\eqref{extend-cond} cannot suffice as a substitute for linearity.

Alternatively, one might consider the ansatz
\begin{equation} \label{new-try}
    \cR \big(a^{}_0 x + a^{}_1 \cR(x) + \ldots
    + a^{}_n \cR^n(x)\big)  \; = \; \cR(x)
\end{equation}
on $K$, where the coefficients $a_i\ge 0$ with $a^{}_0 + \ldots +
a^{}_n=1$ define general convex combinations. With $a_j=1$ (hence $a_i
= 0$ for all $i\neq j$), this comprises $\cR^{j+1}=\cR$ on $K$, hence
also $\cR^2=\cR$, so that this ansatz cannot yield an extension of
\eqref{gen-cond}.

\smallskip
Any meaningful generalization of \eqref{gen-cond} seems to require a
{\em partial\/} linearity of $\cR$, namely on a subset of the convex
combinations used in \eqref{new-try}. To construct one possibility,
fix $n\ge 2$, set $\xi:=\exp(2\pi i/n)$, and define
\begin{equation} \label{def-gfun}
   G^{(n)}_{k} (t) \; := \; \frac{1}{n} \big( e^t + 
   \xi^{k}\, e^{\xi t} + \xi^{2k}\, e^{\xi^2 t} + \ldots
   + \xi^{(n-1)k}\, e^{\xi^{n-1}t}\big)
\end{equation}
for $k\in\ZZ$. Due to $\xi^n=1$, this definition is modulo $n$ in the
lower index $k$. Observing
\begin{equation} \label{sum-of-roots}
   \frac{1}{n}\sum_{m=0}^{n-1} \xi^{m M} \; = \;
   \begin{cases} 1, & \mbox{if $n$ divides $M$} \\ 0, & \mbox{otherwise} 
   \end{cases}
\end{equation}
and expanding the exponential factors in \eqref{def-gfun}, one
finds (for $0\le k < n$)
\begin{eqnarray*}
   G^{(n)}_{k} (t) & = & \frac{1}{n} \sum_{\ell\ge 0}
   \frac{t^\ell}{\ell !} \, \big( 1 + \xi^{k+\ell} +
   \xi^{2(k+\ell)}+ \ldots + \xi^{(n-1)(k+\ell)}\big) \\
   & = & \delta^{}_{k,0} + \sum_{m\ge 1}
   \frac{t^{mn-k}}{(mn-k)!}
\end{eqnarray*}
For $0\le k < n$, these functions have the following elementary properties.
\begin{enumerate}
\item $G^{(n)}_{k} (0) = \delta^{}_{k,0}\,$;
\item $\frac{\rm d}{{\rm d}t} G^{(n)}_{k} (t)
   = G^{(n)}_{k+1} (t)\quad$  (with $G^{(n)}_{n} \equiv G^{(n)}_{0}$);
\item $\sum_{k=0}^{n-1} G^{(n)}_{k} (t) = e^t$.
\end{enumerate}
Moreover, one has
\begin{lemma} \label{asymp}
   Let $n\ge 2$ and $0\le k < n$. Then,
   $\;\lim_{t\to\infty} e^{-t} G^{(n)}_{k} (t) = \frac{1}{n}$.
\end{lemma}
{\sc Proof}:
If $\gamma\in S^1$, one has $\gamma=\cos(\phi) + i \sin(\phi)$ for
some $\phi\in\RR$, hence
\[
   e^{-t} e^{\gamma t} \; = \; e^{(\gamma -1)t} \; = \;
   e^{(\cos(\phi)-1)t} e^{i\sin(\phi)t}.
\]
The absolute value is then $e^{(\cos(\phi)-1)t}$ which tends to $0$ as
$t\to\infty$, unless $\cos(\phi)=1$. Using this term by term in
\eqref{def-gfun} proves the claim. \qed

\smallskip Let us assume, as above, that $\cR(K)\subset K$ and
consider the abstract Cauchy problem \eqref{abstract-Cauchy}. Let us
also assume, for a fixed $n\ge 2$, that $\cR^{n+1}=\cR$ on $K$.
Define
\begin{equation} \label{new-flow}
   \varphi^{}_t \; = \; e^{-t}\, \big(
   \bs{1} + (G^{(n)}_{0} (t) - 1)\cR^n +
   \sum_{k=1}^{n-1} G^{(n)}_{n-k} (t) \cR^k \,\big)
\end{equation}
so that $\varphi^{}_0 = \bs{1}$ and $\varphi^{}_{t}
\xrightarrow{t\to\infty}
\frac{1}{n} (\cR+\cR^2+\ldots+\cR^n)$, as a consequence of
Lemma~\ref{asymp}.

If we now have $[\cR,\varphi^{}_{t}]=0$ for all $t>0$, one can check,
by an explicit calculation, that $\varphi^{}_{t} (x^{}_{0})$ solves
the Cauchy problem \eqref{abstract-Cauchy}.  For $n=2$, it is the
solution of Proposition~\ref{single-rec-group}.

Let us make this approach more concrete for $n=3$. Here, one has 
$\cR^3=\cR$, and the solution $\varphi^{}_{t}(x^{}_{0})$ would read
\[
   x(t) \; = \; e^{-t} x^{}_{0} + e^{-t} \sinh(t)\,\cR(x^{}_{0})
   + e^{-t} (\cosh(t)-1)\,\cR^2(x^{}_{0}).
\]
The coefficients on the right hand side once more admit a
probabilistic interpretation. The term $e^{-t}$ is the probability
that no ``hit'' (by $\cR$) has happened until time $t$, while $e^{-t}
\sinh(t)=\frac{1}{2}(1+e^t)(1-e^{-t})$ (resp.\ $e^{-t} (\cosh(t)-1) =
\frac{1}{2} (1-e^{-t})^2$) is the probability for an odd number of
hits (resp.\ an even number $\ge 2$) until time $t$.

Do such operators $\cR$ exist that are nonlinear? One possibility to
construct an example is the following. For an idempotent $R$ with
$R(K)\subset K$, find a map $\sigma\!:\, R(K)\longrightarrow R(K)$
with $\sigma^n=1$ that commutes with $R$ on $R(K)$ (which is
nontrivial).  If one defines $\cR=\sigma R$, one has $\cR^n=R$ and
$\cR^{n+1}=\cR$ on $K$.

\subsection{Dynamics under compatible recombinators}

Let us go back to the more concrete setting of simple recombination
and look at the situation of compatible semigroups.  If $A\subset L$,
let
\[
   I(A)\; := \; \{\beta\mid \min(A)\le\beta\le\max(A)\}
\]
denote the complete stretch of links associated with $A$. 
\begin{lemma} \label{commute-1}
  Let $A\subset L$ be fixed and let $\alpha\in L$ be given with
  $\alpha\not\in I(A)$. Let $\{\varphi^A_t\mid t\ge 0\}$ denote the
  corresponding positive semigroup according to
  Proposition~$\ref{single-rec-group}$. On $\cM_+(X)$, one then has
   $$\varphi^A_t \circ R^{}_{\alpha} \; = \;
     R^{}_{\alpha} \circ \varphi^A_t $$
   for all $t \ge 0$.
\end{lemma}
{\sc Proof}: In view of the formula for the operators
$\varphi^{A}_{t}$, it is again sufficient to prove the claim on
$\cP(X)$. The extension to $\cM_+(X)$ then follows from the positive
homogeneity of the recombinators.
 
Let $\nu$ be an arbitrary probability measure, and assume that $\alpha
< \min(A)$ (the case $\alpha > \max(A)$ is completely analogous). If
$a\in[0,1]$, one finds
\begin{eqnarray*}
   \lefteqn{R^{}_{\alpha}\big(a\nu + (1-a)R^{}_A(\nu)\big)} \\ 
   & = &
   \big(a\, \pi^{}_{\s\alpha}.\nu + (1-a)\, \pi^{}_{\s\alpha}.R^{}_A(\nu)\big)
   \otimes 
   \big(a\, \pi^{}_{\g\alpha}.\nu + (1-a)\, 
            \pi^{}_{\g\alpha}.R^{}_A(\nu)\big) \\
   & = &
   a^2 R^{}_{\alpha}(\nu) + (1-a)^2 R^{}_{A\cup\{\alpha\}}(\nu) + \\
   & &
   a (1-a) \big( (\pi^{}_{\s\alpha}.R^{}_A(\nu))\otimes (\pi^{}_{\g\alpha}.\nu)
   + (\pi^{}_{\s\alpha}.\nu)\otimes (\pi^{}_{\g\alpha}.R^{}_A(\nu))\big).
\end{eqnarray*}
Observe that, due to our assumptions on $\alpha$, we have
$\pi^{}_{\s\alpha}.\nu = \pi^{}_{\s\alpha}.R^{}_A(\nu)$. Using this
relation twice in the last line above (once in each direction), one
obtains
\begin{eqnarray*}
   \lefteqn{a^2 R^{}_{\alpha}(\nu) + (1-a)^2 R^{}_{A\cup\{\alpha\}}(\nu) +
   a (1-a) \big( R^{}_{\alpha}(\nu) + R^{}_{\alpha}(R^{}_A(\nu))\big) } \\
   & = & a \, R^{}_{\alpha}(\nu) + (1-a)\, R^{}_{A\cup\{\alpha\}}(\nu) 
   \,\; = \;\, \big( a\, \bs{1} + (1-a)\, R^{}_A\big) (R^{}_{\alpha}(\nu))\,.
\end{eqnarray*}
Since $\nu\in\cM_+(X)$ was arbitrary, this implies the claim.  \qed

\begin{coro} \label{commute-2}
   Let $A,B \subset L$ with $I(A)\cap I(B) = \es$. Then, on $\cM_+(X)$,
   $$ R^{}_{B} \circ \varphi^{A}_{t} \; = \; 
      \varphi^{A}_{t} \circ R^{}_{B} \, .$$
\end{coro}
{\sc Proof}: Let $\psi = a\,\bs{1} + (1-a) R^{}_{A}$ for an arbitrary,
but fixed $a\in [0,1]$. Write $B=\{\beta_1,\dots,\beta_r\}$ with
$\beta_1 < \dots < \beta_r$. Then, $R^{}_{B} = \prod_{i=1}^{r}
R^{}_{\beta_i}$ by Proposition~\ref{alg-rules}.4, and a repeated
application of Lemma~\ref{commute-1} gives
%\begin{eqnarray*}
%   R^{}_{B} \circ \psi 
%   & = & \Big(\prod_{i=1}^{r} R^{}_{\beta_i}\Big) \circ \psi 
%   \,\; = \,\; \Big(\prod_{i=1}^{r-1} R^{}_{\beta_i}\Big) \circ
%          \psi \circ R^{}_{\beta_r} \\
%   & = & \dots \,\; = \,\; R^{}_{\beta_1} \circ
%          \psi \circ
%          \Big(\prod_{i=2}^{r} R^{}_{\beta_i}\Big) \\
%   & = &  \psi \circ \Big(\prod_{i=1}^{r} R^{}_{\beta_i}\Big)  
%   \,\; = \,\;  \psi \circ R^{}_{B} 
%\end{eqnarray*}
$$  R^{}_{B} \circ \psi 
    \; = \; \Big(\prod_{i=1}^{r} R^{}_{\beta_i}\Big) \circ \psi 
    \; = \; \dots \; = \;
    \psi \circ \Big(\prod_{i=1}^{r} R^{}_{\beta_i}\Big)  
    \; = \; \psi \circ R^{}_{B}  $$
which is valid on $\cM_+(X)$.
This proves the claim because $\varphi^{A}_{t}$ is of the form 
$( a\bs{1} + (1-a) R^{}_{A})$, with $a\in [0,1]$, for all $t\ge 0$. \qed

\begin{theorem} \label{commute-3}
  Let $A,B \subset L$ with $I(A)\cap I(B) = \es$. Then, the
  corresponding semigroups commute, i.e.,
   $$ \varphi^{A}_{t} \circ \varphi^{B}_{s} \; = \;
      \varphi^{B}_{s} \circ \varphi^{A}_{t} $$
   on $\cM_+(X)$, for all\/ $t,s\ge 0$.
   
   In particular, $\{\varphi^{A}_{t} \circ \varphi^{B}_{s}\mid t,s\ge
   0\}$ defines an Abelian two-parameter semigroup of nonlinear
   positive operators on $\cM_+(X)$.
\end{theorem}
{\sc Proof}: Let $\nu\in\cP(X)$. With Corollary~\ref{commute-2}, one finds
\begin{eqnarray*}
  \lefteqn{\big( (a\bs{1} + (1-a) R^{}_{A}) \circ
           ( b\bs{1} + (1-b) R^{}_{B}) \big) (\nu)} \\
  & = & ab\,\nu + a (1-b) R^{}_{B} (\nu) + (1-a)
        \big( b\bs{1} + (1-b) R^{}_{B}\big) (R^{}_{A}(\nu)) \\
  & = &  ab\,\nu + a (1-b) R^{}_{B} (\nu) + b(1-a) R^{}_{A} (\nu)
         + (1-a) (1-b) R^{}_{A\cup B} (\nu) \, .
\end{eqnarray*}
Since the last expression is symmetric in $(a,A)$ versus $(b,B)$,
the operators in the first line commute. 

Since the operators $\varphi^{A}_{t}$ and $\varphi^{B}_{s}$ are of the
form used in this argument, for all $t,s\ge 0$, the claim is true on
$\cP(X)$. By positive homogeneity of the recombinators, it extends to
all of $\cM_+(X)$.  \qed

\smallskip This allows to formulate our main result, where we write
$\varphi^{L_i}_{t}$ for the semigroup attached to a set $L_i\subset L$
of links.
\begin{theorem}  \label{mixed-semigroups}
   Let $A:=\,\bigcup_{1\le i\le r}\, L_i$ be a subset of $L$, with 
   $I(L_i)\cap I(L_j)=\es$ for all\/ $i\neq j$.
   Then, the semigroups $\{\varphi^{L_i}_{t}\mid t\ge 0\}$ mutually
   commute, and the Cauchy problem
$$ \dot{\omega}\; = \; \sum_{i=1}^{r} \varrho^{}_{L_i}
   \big(R^{}_{L_i} - 1\big) (\omega)\, , $$
with all $\varrho^{}_{L_i} > 0$ and initial condition $\omega^{}_0 \ge 0$, 
has the unique solution
$$
\omega^{}_t \; = \; \Big( \prod_{i=1}^{r} \varphi^{L_i}_{t}\Big) \,
(\omega^{}_0) \, . $$
Asymptotically, $\omega^{}_t
\stackrel{\|.\|\,}{\longrightarrow}R^{}_A(\omega^{}_0)$ as
$t\to\infty$, the convergence, once again, being exponentially fast.
\end{theorem}
{\sc Proof}: Due to commutativity of the participating semigroups by
Theorem~\ref{commute-3}, one can apply Proposition~\ref{single-rec-group}
repeatedly to find
\begin{eqnarray*}
  \dot{\omega}^{}_{t}& = & \Big(\sum_{i=1}^{r}
  \Big( \frac{\rm d}{{\rm d}t}\, \varphi^{L_i}_{t} \Big) \circ
  \prod_{j\neq i} \varphi^{L_j}_{t}\Big) \, (\omega^{}_0) \\
  & = &\Big( \sum_{i=1}^{r}  \varrho^{}_{L_i} \big(R^{}_{L_i} -1\big)
  \circ \prod_{i=1}^{r} \varphi^{L_i}_{t}\Big) \, (\omega^{}_0) \\
  & = & \sum_{i=1}^{r}  \varrho^{}_{L_i} \big(R^{}_{L_i} -1\big)
  (\omega^{}_{t})\, .
\end{eqnarray*}
The convergence result towards the product measure $R^{}_A(\omega^{}_0)$
is a multiple
application of Proposition~\ref{single-rec-group}, together
with Part 4 of Proposition~\ref{alg-rules}. \qed

\begin{coro}  \label{multi-par}
   Under the assumptions of Theorem~$\ref{mixed-semigroups}$, the set 
$$ \{ \prod_{i=1}^{r} \varphi^{L_i}_{t_i} \mid t_i \ge 0\}$$
   forms an Abelian $r$-parameter semigroup of nonlinear
   positive operators on $\cM_+(X)$. Each factor is of the
   form
$$ \varphi^{L_i}_{t_i} \; = \; \exp(-\varrho^{}_{L_i}t^{}_i)\, \bs{1}
   + \big(1 - \exp(-\varrho^{}_{L_i}t^{}_i)\big) R^{}_{L_i} $$
   where $\varrho^{}_{L_i} > 0$ is the intensity of the
   underlying process. \qed
\end{coro}

\smallskip Another obvious consequence is that the forward flow of
Theorem~\ref{mixed-semigroups}, which is a one-parameter semigroup, is
embedded into the $r$-parameter semigroup of
Corollary~\ref{multi-par}.

\subsection{Application to single crossovers and outlook}

In \cite{reco}, the biologically most relevant situation was
investigated where $L$ was written as the disjoint union of all its
elements. This resulted in the ODE
\begin{equation}  \label{old}
   \dot{\omega} \; = \;  \sum_{\alpha\in L} \varrho^{}_{\alpha}
   \big(R^{}_{\alpha} - 1 \big) (\omega)\, .
\end{equation}
By a different approach, it was shown that the corresponding Cauchy
problem with positive initial condition $\omega^{}_0$ has the
solution
\begin{equation} \label{old-sol}
    \omega^{}_t \; = \;  
    \sum_{G\subset L} a^{}_G (t)\, R^{}_{G} (\omega^{}_0)  
\end{equation}
where the coefficient functions are given by
$$  a^{}_G (t) \; = \; \prod_{\alpha\in\overline{G}}
    \exp(-\varrho^{}_{\alpha} t) \prod_{\beta\in G} 
    \big(1-\exp(-\varrho^{}_{\beta} t)\big) .  $$
    
In our new (and more general) formulation, these coefficients can
be seen as the result of expanding the product (over $\alpha\in
L$) of the commuting semigroups $\{\varphi^{\alpha}_{t}\mid t\ge
0\}$.  It is an easy exercise to show that the formula of
Proposition~\ref{single-rec-group} then gives an independent
verification of (\ref{old-sol}).

In this case, a simple combinatorial transformation is possible, namely
\[
   T^{}_G \; := \; \sum_{H\supset G} (-1)^{|H-G|}\, R^{}_H ,
\]
which admits the inverse $R^{}_G := \sum_{H\supset G} T^{}_H$
by M\"obius inversion, compare \cite[Thm.~4.18]{Aigner}. It was shown in
\cite{reco} that the signed measures $T^{}_G (\omega^{}_t)$ satisfy the
{\em linear}\/ ODEs
$$   \frac{\rm d}{{\rm d}t}\, T^{}_G (\omega^{}_t) \; = \;
     - \Big(\sum_{\alpha\in\overline{G}} \varrho^{}_{\alpha}\Big)\,
     T^{}_G (\omega^{}_t) \, . $$
Since $T^{}_G (\omega^{}_t) = b^{}_G (t)\, T^{}_G (\omega^{}_0)$ with
$b^{}_G (t) = \sum_{H\subset G} a^{}_H (t)$, $\omega^{}_t$ of
(\ref{old-sol}) also admits the representation
$$ \omega^{}_t \; = \; 
   \sum_{G\subset L} b^{}_G (t)\, T^{}_{G} (\omega^{}_0) \, . $$
This shows that (\ref{old}) can be transformed to a system
of $2^{|L|}$ linear ODEs (together with a special set of initial
conditions), which, a posteriori, provides an explanation 
for the appearance of the ``almost linear like'' behaviour of the
nonlinear flow of the original equation (\ref{old}). At present, I 
am not aware of any other example of this kind of ``M\"{o}bius
linearization'' in the literature.

\smallskip
A similar observation applies to the situation of 
Theorem~\ref{mixed-semigroups}, which corroborates the
intuition that explicit expressions for nonlinear semigroups
ought to be related to some linear structure. As the above
examples show, there are perhaps more possibilities for such
a connection to be discovered.

Our above analysis revolved around ordered partitions, which make the
treatment rather simple and transparent. It is clear, however, that
there is no principal reason to restrict oneself to this case, and
with some extra effort, similar results might also be possible for
more general partitions.  One difficulty here is to find a good
formulation for the cases where the semigroups (and not just the
recombinators) commute.

\medskip
\subsection*{Acknowledgements}

It is a pleasure to thank Ellen Baake, Uwe Grimm and Manfred Wolff for 
discussions and valuable comments on the manuscript, and the Erwin
Schr\"odinger International Institute in Vienna for hospitality. 
Financial support  from the German academic exchange service (DAAD) 
is gratefully acknowledged.

\end{document}